\RequirePackage{fix-cm}
\RequirePackage{fixltx2e}
\documentclass[oneside,english]{amsart}

\usepackage[latin9]{inputenc}
\usepackage{babel}
\usepackage{url}
\usepackage{amsthm}
\usepackage{amssymb}
\usepackage{esint}
\usepackage[unicode=true,pdfusetitle,
 bookmarks=true,bookmarksnumbered=false,bookmarksopen=false,
 breaklinks=false,pdfborder={0 0 1},backref=false,colorlinks=false]
 {hyperref}
 \makeatletter

\renewcommand{\theequation}{\theequation. \arabic{equation}}
\numberwithin{equation}{section}
\newtheorem{thm}{Theorem}[section]
\newtheorem{cor}{Corollary}[section]

\newtheorem{prop}{Proposition}[section]

\usepackage[numbers,sort&compress]{natbib}

\def\squarebox#1{\hbox to #1{\hfill\vbox to #1{\vfill}}}
\def\qed{\hspace*{\fill}
         \vbox{\hrule\hbox{\vrule\squarebox{.667em}\vrule}\hrule}\smallskip}
\begin{document}\large
\title[A Lauricella hypergeometric series over finite fields]
{\Large A Lauricella hypergeometric series over finite fields}
\author{\small  Bing He}
\address{\small
College of Science, Northwest A\&F University,
   Yangling 712100, Shaanxi, People's Republic of China}
\email{yuhe001@foxmail.com; yuhelingyun@foxmail.com}




\begin{abstract}
\small In this paper we introduce a finite field analogue for a Lauricella hypergeometric series.  An integral  formula for the Lauricella hypergeometric series and its finite field analogue are deduced. Transformation and reduction formulae and  several generating functions for the Lauricella hypergeometric series over finite fields are obtained.  Some of these generalize certain results of Li \emph{et al} and Greene as well as several other known results.
\end{abstract}
\maketitle
\tableofcontents
\section{Introduction}

Let $q$ be a power of a prime and let $\mathbb{F}_{q}$ and $\widehat{\mathbb{F}}^{*}_{q}$  denote the finite field of $q$ elements and  the group of multiplicative characters
of $\mathbb{F}^{*}_{q}$ respectively. Setting $\chi(0)=0$ for all characters, we extend the domain of all characters $\chi$ of $\mathbb{F}^{*}_{q}$ to $\mathbb{F}_{q}.$  Let $\overline{\chi}$ and $\varepsilon$ denote the inverse of $\chi$ and the trivial character respectively. See \cite{BEW} and \cite [Chapter 8]{IR} for more information about characters.

Following \cite{B}, we define the generalized hypergeometric function as
$$ {}_{n+1}F_n \left(\begin{matrix}
a_0, a_1, \ldots , a_{n} \\
b_1, \ldots , b_n \end{matrix}
\bigg| x \right):=\sum_{k=0}^{\infty}\frac{(a_{0})_{k}(a_{1})_{k}\cdots(a_{n})_{k}}{k!(b_{1})_{k}\cdots(b_{n})_{k}}x^{k},$$
where $(z)_{k}$ is the Pochhammer symbol given by
\begin{equation*}
(z)_{0}=1,~(z)_{k}=z(z+1)\cdots(z+k-1)\text{ for } k\geq 1.
\end{equation*}
It was Greene  who in \cite{Gr} developed the theory of hypergeometric functions over finite fields and established numerous transformation and summation identities for hypergeometric series over finite fields which are analogues to those in the classical case. Greene, in particular,  introduced the notation
\begin{equation*}
{}_{2}F_1 \left(\begin{matrix}
A, B \\
C \end{matrix}
\bigg| x \right)^{G}=\varepsilon(x)\frac{BC(-1)}{q}\sum_{y}B(y)\overline{B}C(1-y)\overline{A}(1-xy)
\end{equation*}
for $A,B,C\in \widehat{\mathbb{F}}_{q}$ and $x\in \mathbb{F}_{q},$ that is
a finite field analogue for the integral representation of Gauss hypergeometric series \cite{B}:
\begin{equation*}
{}_{2}F_1 \left(\begin{matrix}
a, b \\
c \end{matrix}
\bigg| x \right)=\frac{\Gamma(c)}{\Gamma(b)\Gamma(c-b)}\int_{0}^{1}t^b(1-t)^{c-b}(1-tx)^{-a}\frac{dt}{t(1-t)},
\end{equation*}
and defined the finite field analogue for the binomial coefficient as
\begin{equation*}
  {A\choose B}^{G}=\frac{B(-1)}{q}J(A,\overline{B}),
\end{equation*}
where $J(\chi,\lambda)$ is the Jacobi sum given by $$J(\chi,\lambda)=\sum_{u}\chi(u)\lambda(1-u).$$ For more details about the finite field analogue for the generalized hypergeometric functions, please see \cite{FL, M, EG}.

In this paper,  for the sake of simplicity, we define the finite field analogue for the binomial coefficient and the classic Gauss hypergeometric series by
\begin{equation*}
  {A\choose B}=q{A\choose B}^{G}=B(-1)J(A,\overline{B}).
\end{equation*}
and
\begin{equation*}
{}_{2}F_1 \left(\begin{matrix}
A, B \\
C \end{matrix}
\bigg| x \right)=q\cdot{}_{2}F_1 \left(\begin{matrix}
A, B \\
C \end{matrix}
\bigg| x \right)^{G}=\varepsilon(x)BC(-1)\sum_{y}B(y)\overline{B}C(1-y)\overline{A}(1-xy),
\end{equation*}
respectively.

There are many interesting double hypergeometric functions in the field of hypergeometric functions. Among these functions, the Appell series $F_{1}$ may be one of the most important functions:
\begin{align*}
F_{1}(a;b,b';c;x,y)&=\sum_{m,n\geq 0}\frac{(a)_{m+n}(b)_{m}(b')_{n}}{m!n!(c)_{m+n}}x^my^n,~|x|<1,~|y|<1.
\end{align*}
See \cite{B, CA, S} for more material about the Appell series.

Inspired by Greene's work, Li \emph{et al} \cite{LLM} gave a finite field analogue for the Appell series $F_{1}$ and
established some transformation and reduction formulas and the generating functions for the function over finite fields.  In that paper, the finite field analogue
for the Appell series $F_{1}$ was given by
\begin{equation*}
F_{1}(A;B,B';C;x,y)=\varepsilon(xy)AC(-1)\sum_{u}A(u)\overline{A}C(1-u)\overline{B}(1-ux)\overline{B'}(1-uy).
\end{equation*}

It was Lauricella \cite{L} who in 1893 generalized the Appell series $F_{1}$ to  the Lauricella hypergeometric series $F^{(n)}_{D}$  that is defined by
\begin{equation*}
F^{(n)}_{D}\left(\begin{matrix}
a; b_{1},\cdots,b_{n} \\
 c \end{matrix}
\bigg| x_{1},\cdots,x_{n}  \right):=\sum_{m_{1}=0}^{\infty}\cdots\sum_{m_{n}=0}^{\infty}\frac{(a)_{m_{1}+\cdots+m_{n}}(b_{1})_{m_{1}}\cdots(b_{n})_{m_{n}}}{(c)_{m_{1}+\cdots+m_{n}}m_{1}!\cdots m_{n}!}x_{1}^{m_{1}}\cdots x_{n}^{m_{n}}.
\end{equation*}
It is clear that
\begin{equation*}
F_{1}(a;b,b';c;x,y)=F^{(2)}_{D}\left(\begin{matrix}
a; b,b^{'} \\
 c \end{matrix}
\bigg| x,y\right)~\text{and}~{}_{2}F_{1}\left(\begin{matrix}
b, a \\
 c \end{matrix}
\bigg| x\right)=F^{(1)}_{D}\left(\begin{matrix}
a; b \\
 c \end{matrix}
\bigg| x\right).
\end{equation*}
so the Lauricella hypergeometric series $F^{(n)}_{D}$ is an $n$-variable extension of the Appell series $F_{1}$ and  the hypergeometric function  ${}_{2}F_{1}.$

Motivated by the work of Greene \cite{Gr} and Li \emph{et al} \cite{LLM}, we give a finite field analogue for  the Lauricella hypergeometric series. Since the Lauricella hypergeometric series $F^{(n)}_{D}$ has a integral representation
\begin{equation*}
F^{(n)}_{D}\left(\begin{matrix}
a; b_{1},\cdots,b_{n} \\
 c \end{matrix}
\bigg| x_{1},\cdots,x_{n}  \right)=\frac{\Gamma(c)}{\Gamma(a)\Gamma(c-a)}\int_{0}^{1}\frac{u^{a-1}(1-u)^{c-a-1}}{(1-x_{1}u)^{b_{1}}\cdots(1-x_{n}u)^{b_{n}}}du,
\end{equation*}
we give the finite field analogue for the Lauricella hypergeometric series  in the following form:
\begin{align*}
  &F^{(n)}_{D}\left(\begin{matrix}
A; B_{1},\cdots,B_{n} \\
 C \end{matrix}
\bigg| x_{1},\cdots,x_{n}\right)\\
&~~~~~~~=\varepsilon(x_{1}\cdots x_{n})AC(-1)\sum_{u}A(u)\overline{A}C(1-u)\overline{B_{1}}(1-x_{1}u)\cdots \overline{B_{n}}(1-x_{n}u),
\end{align*}
where $A,B_{1},\cdots,B_{n},C_{1},\cdots,C_{n}\in \widehat{\mathbb{F}}_{q},~x_{1},\cdots,x_{n}\in \mathbb{F}_{q}$ and  the sum ranges over all the elements of $\mathbb{F}_{q}.$ In the above definition,  the factor $\frac{\Gamma(c)}{\Gamma(a)\Gamma(c-a)}$ is dropped to obtain simpler results. We choose the factor  $\varepsilon(x_{1}\cdots x_{n})AC(-1)$ to get a better expression in terms of binomial coefficients. From the definition of the Lauricella hypergeometric series over finite fields, we know that
\begin{equation*}
  F_{1}(A;B,B';C;x,y)=F^{(2)}_{D}\left(\begin{matrix}
A; B,B' \\
 C \end{matrix}
\bigg| x,y\right)~\text{and}~{}_{2}F_1 \left(\begin{matrix}
A, B \\
C \end{matrix}
\bigg| x \right)=F^{(1)}_{D}\left(\begin{matrix}
B; A \\
 C \end{matrix}
\bigg| x\right).
\end{equation*}
Then the Lauricella hypergeometric series over finite fields can be regarded as an $n$-variable extension of the finite field analogues for the Appell series $F_{1}$ and  the hypergeometric function  ${}_{2}F_{1}.$

From Theorem \ref{t1} or the definition, we know that the Lauricella hypergeometric series  over finite fields
\begin{equation*}
  F^{(n)}_{D}\left(\begin{matrix}
A; B_{1},\cdots,B_{n} \\
 C \end{matrix}
\bigg| x_{1},\cdots,x_{n}\right)
\end{equation*}
is invariant under permutation of the subscripts $1,2,\cdots, n,$ namely, it is invariant under permutation of the $B'$s and $x'$s together.

The aim of this paper is to give several transformation and reduction formulas and  the generating functions for the Lauricella hypergeometric series over finite fields. We know that the Lauricella hypergeometric series over finite fields is an $n$-variable extension of the finite field analogues for the Appell series $F_{1}$ and  the hypergeometric function  ${}_{2}F_{1}.$ So some of the results in this paper are generalizations of certain results in \cite{LLM, Gr} and some other known results. For example,   \cite [Theorem 1.3]{LLM} and \cite [Theorem 3.6]{Gr} are special cases of Theorem \ref{t1}.

We will give another expression for the Lauricella hypergeometric series over finite fields in the next section.  An integral  formula for the Lauricella hypergeometric series and its finite field analogue are deduced in Section 3. In Section 4, several transformation and reduction formulae for the Lauricella hypergeometric series  over finite fields will be given.  The last section is devoted to some generating functions for the Lauricella hypergeometric series over finite fields.
\section{Another expression}
In this section we give another expression for the Lauricella hypergeometric series over finite fields.
\begin{thm}\label{t1}For $A,B_{1},\cdots,B_{n},C\in \widehat{\mathbb{F}}_{q}$ and $x_{1},\cdots,x_{n}\in \mathbb{F}_{q},$ we have
\begin{align*}
  &F^{(n)}_{D}\left(\begin{matrix}
A; B_{1},\cdots,B_{n} \\
 C \end{matrix}
\bigg| x_{1},\cdots,x_{n}\right)\\
&~~=\frac{1}{(q-1)^n}\sum_{\chi_{1},\cdots, \chi_{n}}{A\chi_{1}\cdots\chi_{n}\choose C\chi_{1}\cdots\chi_{n}}{B_{1}\chi_{1}\choose \chi_{1}}\cdots {B_{n}\chi_{n}\choose \chi_{n}}\chi_{1}(x_{1})\cdots \chi_{n}(x_{n}),
\end{align*}
where each sum ranges over all multiplicative characters of $\mathbb{F}_{q}.$
\end{thm}
To carry out our study, we need some auxiliary results which will be used  in the sequel.

The results in the  following proposition  follows readily from some properties of Jacobi sums.
\begin{prop} \emph{(See \cite [(2.6), (2.8) and (2.12)]{Gr})} If $A,B\in \widehat{\mathbb{F}}_{q},$ then
\begin{align}
  {A\choose B}&={A\choose A\overline{B}},\label{f2}\\
  {A\choose B}&={\overline{B}\choose \overline{A}}AB(-1),\label{f3}\\
{A\choose \varepsilon}&={A\choose A}=-1+(q-1)\delta(A),\label{f4}
\end{align}
where $\delta(\chi)$ is a function on characters given by
\begin{equation*}
\delta(\chi)=\left\{
               \begin{array}{ll}
                 1  & \hbox{if $\chi=\varepsilon$} \\
                 0  & \hbox{otherwise}
               \end{array}.
             \right.
\end{equation*}
\end{prop}
The following result is a finite field analogue for the well-known identity $${a\choose b}{c\choose a}={c\choose b}{c-b\choose a-b}.$$
\begin{prop}\label{p2-2}\emph{(See \cite [(2.15)]{Gr})}For $A,B,C\in \widehat{\mathbb{F}}_{q},$ we have
\begin{equation*}
  {A\choose B}{C\choose A}={C\choose B}{C\overline{B}\choose A\overline{B}}-(q-1)B(-1)\delta(A)+(q-1)AB(-1)\delta(B\overline{C}).
\end{equation*}
\end{prop}
The following result is also very important in the derivation of Theorem \ref{t1}.
\begin{thm}\emph{(Binomial theorem over finite fields, see \cite [(2.10)]{Gr})} For $A\in \widehat{\mathbb{F}}_{q}$ and $x\in \mathbb{F}_{q},$ we have
\begin{equation*}
  \overline{A}(1-x)=\delta(x)+\frac{1}{q-1}\sum_{\chi}{A\chi\choose \chi}\chi(x),
\end{equation*}
where the sum ranges over all multiplicative characters of $\mathbb{F}_{q}$ and $\delta(x)$ is a function on $\mathbb{F}_{q}$ given by
\begin{equation*}
\delta(x)=\left\{
            \begin{array}{ll}
              1 & \hbox{if $x=0$} \\
              0 & \hbox{if $x\neq 0$}
            \end{array}.
          \right.
\end{equation*}
\end{thm}

We are now ready to prove Theorem \ref{t1}.\\
\emph{Proof of Theorem \ref{t1}.} From the binomial theorem over finite fields, we know that for $1\leq j\leq n,$
\begin{equation*}
\overline{B_{j}}(1-x_{j}u)=\delta(x_{j}u)+\frac{1}{q-1}\sum_{\chi_{j}}{B_{j}\chi_{j}\choose \chi_{j}}\chi_{j}(x_{j}u).
\end{equation*}
Then, by the fact that $\varepsilon(x_{j})\delta(x_{j}u)A(u)=0$ for $1\leq j\leq n,$
\begin{align*}
  &F^{(n)}_{D}\left(\begin{matrix}
A; B_{1},\cdots,B_{n} \\
 C \end{matrix}
\bigg| x_{1},\cdots,x_{n}\right)\\
&=\varepsilon(x_{1}\cdots x_{n})AC(-1)\sum_{u}A(u)\overline{A}C(1-u)\\
&\cdot \left(\delta(x_{1}u)+\frac{1}{q-1}\sum_{\chi_{1}}{B_{1}\chi_{1}\choose \chi_{1}}\chi_{1}(x_{1}u)\right)\cdots \left(\delta(x_{n}u)+\frac{1}{q-1}\sum_{\chi_{n}}{B_{n}\chi_{n}\choose \chi_{n}}\chi_{n}(x_{n}u)\right)\\
&=\frac{AC(-1)}{(q-1)^n}\sum_{\chi_{1},\cdots,\chi_{n}}{B_{1}\chi_{1}\choose \chi_{1}}\cdots{B_{n}\chi_{n}\choose \chi_{n}}\chi_{1}(x_{1})\cdots \chi_{n}(x_{n})\sum_{u}A\chi_{1}\cdots\chi_{n}(u)\overline{A}C(1-u)\\
&=\frac{1}{(q-1)^n}\sum_{\chi_{1},\cdots,\chi_{n}}{A\chi_{1}\cdots\chi_{n}\choose A\overline{C}}{B_{1}\chi_{1}\choose \chi_{1}}\cdots{B_{n}\chi_{n}\choose \chi_{n}}\chi_{1}(x_{1})\cdots \chi_{n}(x_{n}),
\end{align*}
which, by \eqref{f2}, implies that
\begin{align*}
 &F^{(n)}_{D}\left(\begin{matrix}
A; B_{1},\cdots,B_{n} \\
 C \end{matrix}
\bigg| x_{1},\cdots,x_{n}\right)\\
&~~=\frac{1}{(q-1)^n}\sum_{\chi_{1},\cdots,\chi_{n}}{A\chi_{1}\cdots\chi_{n}\choose C\chi_{1}\cdots\chi_{n}}{B_{1}\chi_{1}\choose \chi_{1}}\cdots{B_{n}\chi_{n}\choose \chi_{n}}\chi_{1}(x_{1})\cdots \chi_{n}(x_{n}).
\end{align*}
This completes the proof of Theorem \ref{t1}.\qed
\section{An integral  formula and its finite field analogue}
In this section we derive an integral  formula for the Lauricella hypergeometric series $F^{(n)}_{D},$ which relating $F^{(n)}_{D}$ to $F^{(n-1)}_{D},$ and then give its finite field analogue. In addition, a finite field analogue for a summation formula on the Lauricella hypergeometric series $F^{(n)}_{D}$ is also deduced.

The following theorem  gives a result which is analogous to \cite[(3.12)]{Gr}.
\begin{thm}\label{t5-1}If $a, b_{1},\cdots,b_{n},c$ and $x_{1},\cdots,x_{n}$ are complex numbers with $Re(b_{1})>0$ and $Re(b_{2})>0,$ then
\begin{align*}
&B(b_{1},b_{2})F^{(n)}_{D}\left(\begin{matrix}
a; b_{1},\cdots,b_{n} \\
 c \end{matrix}
\bigg| x_{1},\cdots,x_{n}  \right)\\
&=\int_{0}^{1}u^{b_{1}-1}(1-u)^{b_{2}-1}F^{(n-1)}_{D}\left(\begin{matrix}
a; b_{1}+b_{2},b_{3},\cdots,b_{n} \\
 c \end{matrix}
\bigg| ux_{1}+(1-u)x_{2},x_{3},\cdots,x_{n}  \right)du,
\end{align*}
where $B(x,y)$ is the beta integral given for $Re(x)>0,Re(y)>0$ by  \cite{An}
\begin{equation*}
B(x,y)=\int_{0}^{1}t^{x-1}(1-t)^{y-1}dt.
\end{equation*}
\end{thm}
\noindent{\it Proof.} It is easily seen from the binomial theorem that
\begin{equation*}
(ux_{1}+(1-u)x_{2})^{m}=\sum_{0\leq m_{1}\leq m}{m\choose m_{1}}u^{m_{1}}(1-u)^{m-m_{1}}x_{1}^{m_{1}}x_{2}^{m-m_{1}}.
\end{equation*}
Then, by \cite [Thoeorem 1.1.4]{An}
\begin{align*}
  &\int_{0}^{1}u^{b_{1}-1}(1-u)^{b_{2}-1}(ux_{1}+(1-u)x_{2})^{m}du\\
  &=\sum_{0\leq m_{1}\leq m}{m\choose m_{1}}x_{1}^{m_{1}}x_{2}^{m-m_{1}}\int_{0}^{1}u^{m_{1}+b_{1}-1}(1-u)^{m-m_{1}+b_{2}-1}du\\
  &=\sum_{0\leq m_{1}\leq m}{m\choose
  m_{1}}x_{1}^{m_{1}}x_{2}^{m-m_{1}}B(m_{1}+b_{1},m-m_{1}+b_{2})\\
  &=\frac{m!}{(b_{1}+b_{2})_{m}}B(b_{1},b_{2})\sum_{0\leq m_{1}\leq m}\frac{(b_{1})_{m_{1}}(b_{2})_{m-m_{1}}}{m_{1}!(m-m_{1})!}x_{1}^{m_{1}}x_{2}^{m-m_{1}}.
\end{align*}
Using the above  identity in the integral at the right side and making the substitution $m-m_{1}\rightarrow m_{2},$ we obtain
\begin{align*}
  &\int_{0}^{1}u^{b_{1}-1}(1-u)^{b_{2}-1}F^{(n-1)}_{D}\left(\begin{matrix}
a; b_{1}+b_{2},b_{3},\cdots,b_{n} \\
 c \end{matrix}
\bigg| ux_{1}+(1-u)x_{2},x_{3},\cdots,x_{n}  \right)du\\
&=\sum_{m=0}^{\infty}\sum_{m_{3}=0}^{\infty}\cdots
\sum_{m_{n}=0}^{\infty}\frac{(a)_{m+m_{3}+\cdots+m_{n}}(b_{1}+b_{2})_{m}(b_{3})_{m_{3}}\cdots(b_{n})_{m_{n}}}{(c)_{m+m_{3}+\cdots+m_{n}}m!m_{3}!\cdots m_{n}!}x_{3}^{m_{3}}\cdots x_{n}^{m_{n}}\\
&~\cdot \int_{0}^{1}u^{b_{1}-1}(1-u)^{b_{2}-1}(ux_{1}+(1-u)x_{2})^{m}du\\
&=B(b_{1},b_{2})\sum_{m_{1}=0}^{\infty}\cdots
\sum_{m_{n}=0}^{\infty}\frac{(a)_{m_{1}+\cdots+m_{n}}(b_{1})_{m}\cdots(b_{n})_{m_{n}}}{(c)_{m_{1}+\cdots+m_{n}}m_{1}!\cdots m_{n}!}x_{1}^{m_{1}}\cdots x_{n}^{m_{n}}\\
&=B(b_{1},b_{2})F^{(n)}_{D}\left(\begin{matrix}
a; b_{1},\cdots,b_{n} \\
 c \end{matrix}
\bigg| x_{1},\cdots,x_{n}  \right).
\end{align*}
We comlete the proof of Theorem  \ref{t5-1}.\qed

Taking $n=2$ in Theorem  \ref{t5-1},  we have the following result.
\begin{cor}If $a, b_{1},b_{2},c$ and $x_{1},x_{2}$ are complex numbers with $Re(b_{1})>0$ and $Re(b_{2})>0,$ then
\begin{align*}
&B(b_{1},b_{2})F_{1}(a; b_{1},b_{2};c ;x_{1},x_{2})=\int_{0}^{1}u^{b_{1}-1}(1-u)^{b_{2}-1}{}_{2}F_{1}\left(\begin{matrix}
 a, b_{1}+b_{2}\\
 c \end{matrix}
\bigg| ux_{1}+(1-u)x_{2}\right)du.
\end{align*}
\end{cor}

When $c=b_{1}+\cdots+b_{n},$ Theorem \ref{t5-1} reduces to \cite[(7.8)]{C}.

We now give the finite field analogue for Theorem \ref{t5-1}.
\begin{thm}\label{t5-4}For $A,B_{3},\cdots,B_{n},C\in \widehat{\mathbb{F}}_{q}, B_{1}, B_{2}\in \widehat{\mathbb{F}}_{q}\backslash \{\varepsilon\}$ and $x_{1},\cdots,x_{n}\in \mathbb{F}_{q},$ we have
\begin{align*}
  &\varepsilon(x_{1}x_{2})\sum_{u}B_{1}(u)B_{2}(1-u)F^{(n-1)}_{D}\left(\begin{matrix}
A; B_{1}B_{2},B_{3},\cdots,B_{n} \\
 C \end{matrix}
\bigg| ux_{1}+(1-u)x_{2},x_{3,}\cdots,x_{n}\right)\\
&={\overline{B_{1}}\overline{B_{2}} \choose \overline{B_{1}}}F^{(n)}_{D}\left(\begin{matrix}
A; B_{1},\cdots,B_{n} \\
 C \end{matrix}
\bigg| x_{1},\cdots,x_{n}\right)\\
&~-\varepsilon(x_{1}x_{2})B_{1}(-1)\overline{B_{1}}\overline{B_{2}}(x_{1}-x_{2})F^{(n-2)}_{D}\left(\begin{matrix}
A\overline{B_{1}}\overline{B_{2}}; B_{3},\cdots,B_{n} \\
 C\overline{B_{1}}\overline{B_{2}} \end{matrix}
\bigg| x_{3},\cdots,x_{n}\right)\\
&~-B_{1}(x_{2})B_{2}(-x_{1})\overline{B_{1}}\overline{B_{2}}(x_{2}-x_{1})F^{(n-2)}_{D}\left(\begin{matrix}
A; B_{3},\cdots,B_{n} \\
 C \end{matrix}
\bigg| x_{3},\cdots,x_{n}\right).
\end{align*}
\end{thm}
\noindent{\it Proof.} It is easily known from the binomial theorem over finite fields that for $u,x_{1}\in \mathbb{F}^{*}_{q},$ we have
\begin{equation*}
  \chi(ux_{1}+(1-u)x_{2})=\frac{1}{q-1}\sum_{\chi_{1}}{\chi\choose \chi_{1}}\chi_{1}(ux_{1})\chi\overline{\chi_{1}}((1-u)x_{2}).
\end{equation*}
Then
\begin{align*}
  &\varepsilon(x_{1})\sum_{u}B_{1}(u)B_{2}(1-u)\chi(ux_{1}+(1-u)x_{2})\\
  &=\frac{1}{q-1}\sum_{\chi_{1}}{\chi\choose \chi_{1}}\chi_{1}(x_{1})\chi\overline{\chi_{1}}(x_{2})\sum_{u}B_{1}\chi_{1}(u)B_{2}\chi\overline{\chi_{1}}(1-u)\\
  &=\frac{1}{q-1}\sum_{\chi_{1}}{\chi\choose \chi_{1}}\chi_{1}(x_{1})\chi\overline{\chi_{1}}(x_{2})B_{2}\chi\chi_{1}(-1){B_{1}\chi_{1}\choose \overline{B_{2}}\overline{\chi}\chi_{1}}.
\end{align*}
We use the above identity in the summation $\sum_{u}$ and replace $\chi\overline{\chi_{1}}$ by $\chi_{2}$ to find
\begin{align}
&\varepsilon(x_{1}x_{2})\sum_{u}B_{1}(u)B_{2}(1-u)F^{(n-1)}_{D}\left(\begin{matrix}
A; B_{1}B_{2},B_{3},\cdots,B_{n} \\
 C \end{matrix}
\bigg| ux_{1}+(1-u)x_{2},x_{3,}\cdots,x_{n}\right)\label{pf54-1}\\
&=\frac{\varepsilon(x_{1}x_{2})}{(q-1)^{n-1}}\sum_{\chi,\chi_{3},\cdots, \chi_{n}}{A\chi\chi_{3}\cdots\chi_{n}\choose C\chi \chi_{3}\cdots\chi_{n}}{B_{1}B_{2}\chi\choose \chi}{B_{3}\chi_{3}\choose \chi_{3}}\cdots {B_{n}\chi_{n}\choose \chi_{n}}\chi_{3}(x_{3})\cdots \chi_{n}(x_{n})\notag\\
&~\cdot \sum_{u}B_{1}(u)B_{2}(1-u)\chi(ux_{1}+(1-u)x_{2})\notag\\
&=\frac{1}{(q-1)^{n}}\sum_{\chi_{1},\cdots, \chi_{n}}B_{2}\chi_{2}(-1){B_{1}\chi_{1}\choose \overline{B_{2}}\overline{\chi_{2}}}{B_{1}B_{2}\chi_{1}\chi_{2}\choose \chi_{1}\chi_{2}}{\chi_{1}\chi_{2}\choose \chi_{1}}\notag\\
&~\cdot {A\chi_{1}\cdots\chi_{n}\choose C\chi_{1}\cdots\chi_{n}}{B_{3}\chi_{3}\choose \chi_{3}}\cdots {B_{n}\chi_{n}\choose \chi_{n}}\chi_{1}(x_{1})\cdots \chi_{n}(x_{n}).\notag
\end{align}
It follows from Proposition \ref{p2-2} that
 \begin{align*}
   {B_{1}B_{2}\chi_{1}\chi_{2}\choose \chi_{1}\chi_{2}}{\chi_{1}\chi_{2}\choose \chi_{1}}&={B_{1}B_{2}\chi_{1}\chi_{2}\choose \chi_{1}}{B_{1}B_{2}\chi_{2}\choose \chi_{2}}-(q-1)\chi_{1}(-1)\delta(\chi_{1}\chi_{2})\\
   &~+(q-1)\chi_{2}(-1)\delta(B_{1}B_{2}\chi_{2}),
 \end{align*}
 \begin{align*}
 {B_{1}\chi_{1}\choose B_{1}B_{2}\chi_{1}\chi_{2}}{B_{1}B_{2}\chi_{1}\chi_{2}\choose \chi_{1}}&={B_{1}\chi_{1}\choose \chi_{1}}{B_{1}\choose B_{1}B_{2}\chi_{2}}-(q-1)\chi_{1}(-1)\delta(B_{1}B_{2}\chi_{1}\chi_{2})
 \end{align*}
 and
 \begin{align*}
 {B_{2}\chi_{2}\choose B_{1}B_{2}\chi_{2}}{B_{1}B_{2}\chi_{2}\choose \chi_{2}}&={B_{2}\chi_{2}\choose \chi_{2}}{B_{2}\choose B_{1}B_{2}}-(q-1)\chi_{2}(-1)\delta(B_{1}B_{2}\chi_{2}).
 \end{align*}
 Then, by \eqref{f2} and \eqref{f3}
 \begin{align*}
&{B_{1}\chi_{1}\choose \overline{B_{2}}\overline{\chi_{2}}}{B_{1}B_{2}\chi_{1}\chi_{2}\choose \chi_{1}\chi_{2}}{\chi_{1}\chi_{2}\choose \chi_{1}}={B_{1}\chi_{1}\choose B_{1}B_{2}\chi_{1}\chi_{2}}{B_{1}B_{2}\chi_{1}\chi_{2}\choose \chi_{1}}{B_{1}B_{2}\chi_{2}\choose \chi_{2}}\\
&~-(q-1)\chi_{1}(-1)\delta(\chi_{1}\chi_{2}){B_{1}\chi_{1}\choose B_{1}B_{2}\chi_{1}\chi_{2}}+(q-1)\chi_{2}(-1)\delta(B_{1}B_{2}\chi_{2}){B_{1}\chi_{1}\choose B_{1}B_{2}\chi_{1}\chi_{2}}\\
&=B_{1}B_{2}\chi_{2}(-1){B_{1}\chi_{1}\choose \chi_{1}}{B_{2}\chi_{2}\choose B_{1}B_{2}\chi_{2}}{B_{1}B_{2}\chi_{2}\choose \chi_{2}}-(q-1)\chi_{1}(-1)\delta(B_{1}B_{2}\chi_{1}\chi_{2}){B_{1}B_{2}\chi_{2}\choose \chi_{2}}\\
&~-(q-1)\chi_{1}(-1)\delta(\chi_{1}\chi_{2}){B_{1}\chi_{1}\choose B_{1}B_{2}\chi_{1}\chi_{2}}+(q-1)\chi_{2}(-1)\delta(B_{1}B_{2}\chi_{2}){B_{1}\chi_{1}\choose B_{1}B_{2}\chi_{1}\chi_{2}}\\
&=B_{2}\chi_{2}(-1){B_{1}\chi_{1}\choose \chi_{1}}{B_{2}\chi_{2}\choose \chi_{2}}{\overline{B_{1}}\overline{B_{2}} \choose \overline{B_{1}}}-(q-1)\chi_{1}(-1)\delta(B_{1}B_{2}\chi_{1}\chi_{2}){B_{1}B_{2}\chi_{2}\choose \chi_{2}}\\
&~-(q-1)\chi_{1}(-1)\delta(\chi_{1}\chi_{2}){B_{1}\chi_{1}\choose B_{1}B_{2}\chi_{1}\chi_{2}},
 \end{align*}
 where in the last step we have cancelled two terms
 \begin{equation*}
   -(q-1)\delta(B_{1}B_{2}\chi_{2})B_{1}B_{2}(-1){B_{1}\chi_{1}\choose \chi_{1}}~\text{and}~(q-1)\chi_{2}(-1)\delta(B_{1}B_{2}\chi_{2}){B_{1}\chi_{1}\choose B_{1}B_{2}\chi_{1}\chi_{2}}.
 \end{equation*}
 Applying the above identity in \eqref{pf54-1}, we obtain
 \begin{align}
 &\varepsilon(x_{1}x_{2})\sum_{u}B_{1}(u)B_{2}(1-u)F^{(n-1)}_{D}\left(\begin{matrix}
A; B_{1}B_{2},B_{3},\cdots,B_{n} \\
 C \end{matrix}
\bigg| ux_{1}+(1-u)x_{2},x_{3,}\cdots,x_{n}\right)\label{pf54-2}\\
&={\overline{B_{1}}\overline{B_{2}} \choose \overline{B_{1}}}\frac{1}{(q-1)^n}\sum_{\chi_{1},\cdots, \chi_{n}}{A\chi_{1}\cdots\chi_{n}\choose C\chi_{1}\cdots\chi_{n}}{B_{1}\chi_{1}\choose \chi_{1}}\cdots {B_{n}\chi_{n}\choose \chi_{n}}\chi_{1}(x_{1})\cdots \chi_{n}(x_{n})\notag\\
&~-\frac{\overline{B_{1}}(-x_{1})\overline{B_{2}}(x_{1})}{(q-1)^n}\sum_{\chi_{3},\cdots, \chi_{n}}{A\overline{B_{1}}\overline{B_{2}}\chi_{3}\cdots\chi_{n}\choose C\overline{B_{1}}\overline{B_{2}}\chi_{3}\cdots\chi_{n}}{B_{3}\chi_{3}\choose \chi_{3}}\cdots {B_{n}\chi_{n}\choose \chi_{n}}\chi_{3}(x_{3})\cdots \chi_{n}(x_{n})\notag\\
&~\cdot \sum_{\chi_{2}}{B_{1}B_{2}\chi_{2}\choose \chi_{2}}\chi_{2}\left(\frac{x_{2}}{x_{1}}\right)\notag\\
&~-\frac{B_{2}(-1)\varepsilon(x_{2})}{(q-1)^{n-1}}\sum_{\chi_{3},\cdots, \chi_{n}}{A\chi_{3}\cdots\chi_{n}\choose C\chi_{3}\cdots\chi_{n}}{B_{3}\chi_{3}\choose \chi_{3}}\cdots {B_{n}\chi_{n}\choose \chi_{n}}\chi_{3}(x_{3})\cdots \chi_{n}(x_{n})\notag\\
&~\cdot \sum_{\chi_{1}}{B_{1}\chi_{1}\choose \overline{B_{2}}\chi_{1}}\chi_{1}\left(\frac{x_{1}}{x_{2}}\right).\notag
 \end{align}
 From \cite [(2.11)]{Gr} we know that for any $A,~B\in \widehat{\mathbb{F}}_{q}$ and $x\in \mathbb{F}_{q},$
\begin{equation}\label{pf31-1}
\sum_{\chi}{A\chi\choose B\chi}\chi(x)=(q-1)\overline{B}(x)\overline{A}B(1-x).
\end{equation}
Using \eqref{pf31-1} in  \eqref{pf54-2} and simplifying yields
\begin{align*}
  &\varepsilon(x_{1}x_{2})\sum_{u}B_{1}(u)B_{2}(1-u)F^{(n-1)}_{D}\left(\begin{matrix}
A; B_{1}B_{2},B_{3},\cdots,B_{n} \\
 C \end{matrix}
\bigg| ux_{1}+(1-u)x_{2},x_{3,}\cdots,x_{n}\right)\\
&={\overline{B_{1}}\overline{B_{2}} \choose \overline{B_{1}}}F^{(n)}_{D}\left(\begin{matrix}
A; B_{1},\cdots,B_{n} \\
 C \end{matrix}
\bigg| x_{1},\cdots,x_{n}\right)\\
&~-\varepsilon(x_{1}x_{2})B_{1}(-1)\overline{B_{1}}\overline{B_{2}}(x_{1}-x_{2})F^{(n-2)}_{D}\left(\begin{matrix}
A\overline{B_{1}}\overline{B_{2}}; B_{3},\cdots,B_{n} \\
 C\overline{B_{1}}\overline{B_{2}} \end{matrix}
\bigg| x_{3},\cdots,x_{n}\right)\\
&~-B_{1}(x_{2})B_{2}(-x_{1})\overline{B_{1}}\overline{B_{2}}(x_{2}-x_{1})F^{(n-2)}_{D}\left(\begin{matrix}
A; B_{3},\cdots,B_{n} \\
 C \end{matrix}
\bigg| x_{3},\cdots,x_{n}\right).
\end{align*}
This concludes the proof of Theorem \ref{t5-4}.\qed

Putting $n=2$ in Theorem \ref{t5-4}, we arrive at
\begin{cor}For $A,C\in \widehat{\mathbb{F}}_{q}, B,B'\in \widehat{\mathbb{F}}_{q}\backslash \{\varepsilon\}$ and $x,y\in \mathbb{F}_{q},$ we have
\begin{align*}
  &\varepsilon(xy)\sum_{u}B(u)B'(1-u){}_{2}F_{1}\left(\begin{matrix}
BB',A\\
 C \end{matrix}
\bigg| ux+(1-u)y\right)\\
&={\overline{B}\overline{B'} \choose \overline{B}}
F_{1}(A; B,B';C; x,y)-\varepsilon(xy)B(-1)\overline{B}\overline{B'}(x-y)\\
&~-B(y)B'(-x)\overline{B}\overline{B'}(y-x).
\end{align*}
\end{cor}

Using the binomial theorem in the integral representation for the Lauricella hypergeometric series and then simplifying, we can get a summation formula connecting  $F^{(n)}_{D}$ and  $F^{(n-1)}_{D}.$
\begin{thm}If $a, b_{1},\cdots,b_{n},c$ and $x_{1},\cdots,x_{n}$ are complex numbers with $Re(a)>0$ and $Re(c-a)>0,$ then
\begin{align*}
 F^{(n)}_{D}\left(\begin{matrix}
a; b_{1},\cdots,b_{n} \\
 c \end{matrix}
\bigg| x_{1},\cdots,x_{n}  \right)=\sum_{k=0}^{\infty}\frac{(a)_{k}(b_{n})_{k}}{k!(c)_{k}}x_{n}^{k}F^{(n-1)}_{D}\left(\begin{matrix}
a+k; b_{1},\cdots,b_{n-1} \\
 c+k \end{matrix}
\bigg| x_{1},\cdots,x_{n-1}  \right).
\end{align*}
\end{thm}
A finite field analogue for the above theorem also holds.
\begin{thm}\label{t5-2}For $A,B_{1},\cdots,B_{n},C\in \widehat{\mathbb{F}}_{q}$ and $x_{1},\cdots,x_{n}\in \mathbb{F}_{q},$ we have
\begin{align*}
  &F^{(n)}_{D}\left(\begin{matrix}
A; B_{1},\cdots,B_{n} \\
 C \end{matrix}
\bigg| x_{1},\cdots,x_{n}\right)\\
&~~~~~=\frac{1}{q-1}\sum_{\chi} {B_{n}\chi\choose \chi}\chi(x_{n})F^{(n-1)}_{D}\left(\begin{matrix}
A\chi; B_{1},\cdots,B_{n-1} \\
 C\chi \end{matrix}
\bigg| x_{1},\cdots,x_{n-1}\right).
\end{align*}
\end{thm}
\noindent{\it Proof.} It follows from the binomial theorem over finite fields that
\begin{equation*}
  \overline{B_{n}}(1-x_{n}u)=\delta(x_{n}u)+\frac{1}{q-1}\sum_{\chi} {B_{n}\chi\choose \chi}\chi(x_{n}u).
\end{equation*}
Using the above identity in the definition of the Lauricella hypergeometric series over finite fields and by the fact that $\varepsilon(x_{n})A(u)\delta(x_{n}u)=0,$ we have
\begin{align*}
  &F^{(n)}_{D}\left(\begin{matrix}
A; B_{1},\cdots,B_{n} \\
 C \end{matrix}
\bigg| x_{1},\cdots,x_{n}\right)\\
&=\frac{\varepsilon(x_{1}\cdots x_{n-1})AC(-1)}{q-1}\sum_{\chi} {B_{n}\chi\choose \chi}\chi(x_{n}u)\\
&~~\cdot\sum_{u}A(u)\overline{A}C(1-u)\overline{B_{1}}(1-x_{1}u)\cdots \overline{B_{n-1}}(1-x_{n-1}u)\\
&=\frac{1}{q-1}\sum_{\chi} {B_{n}\chi\choose \chi}\chi(x_{n}) \varepsilon(x_{1}\cdots x_{n-1})AC(-1)\\
&~\cdot \sum_{u}A\chi(u)\overline{A}C(1-u)\overline{B_{1}}(1-x_{1}u)\cdots \overline{B_{n-1}}(1-x_{n-1}u),
\end{align*}
from which the result follows. This finishes the proof of Theorem \ref{t5-2}.\qed

We set $n=2$ in Theorem \ref{t5-2} to get
\begin{cor}\label{c5-3}For $A,B,B',C\in \widehat{\mathbb{F}}_{q}$ and $x,y\in \mathbb{F}_{q},$ we have
\begin{align*}
F_{1}(A; B,B';C; x,y)=\frac{1}{q-1}\sum_{\chi} {B'\chi\choose \chi}\chi(y){}_{2}F_{1}\left(\begin{matrix}
B,A\chi\\
C\chi \end{matrix}
\bigg|x\right).
\end{align*}
\end{cor}
\cite [Corollary 1.1, (1.5)]{LLM} is a special case of Corollary \ref{c5-3}.
\section{Reduction and Transformation formulae}
In this section we give some reduction and transformation formulae for the Lauricella hypergeometric series  over finite fields.

From the definition of the Lauricella hypergeometric series $F^{(n)}_{D},$ we know that
\begin{align*}
  F^{(n)}_{D}\left(\begin{matrix}
a; b_{1},\cdots,b_{n-1}, 0 \\
 c \end{matrix}
\bigg| x_{1},\cdots,x_{n}  \right)=F^{(n-1)}_{D}\left(\begin{matrix}
a; b_{1},\cdots,b_{n-1} \\
 c \end{matrix}
\bigg| x_{1},\cdots,x_{n-1}  \right).
\end{align*}
We now give a finite field analogue for the above identity.
\begin{thm}\label{t31}For $A,B_{1},\cdots,B_{n-1},C\in \widehat{\mathbb{F}}_{q}$ and $x_{1},\cdots,x_{n}\in \mathbb{F}_{q},$ we have
\begin{align*}
  &F^{(n)}_{D}\left(\begin{matrix}
A; B_{1},\cdots,B_{n-1}, \varepsilon \\
 C \end{matrix}
\bigg| x_{1},\cdots,x_{n}\right)=\varepsilon(x_{n})F^{(n-1)}_{D}\left(\begin{matrix}
A; B_{1},\cdots,B_{n-1} \\
 C \end{matrix}
\bigg| x_{1},\cdots,x_{n-1}\right)\\
&~-\varepsilon (x_{1}\cdots x_{n-1})B_{1}\cdots B_{n-1}\overline{C}(x_{n})\overline{A}C(1-x_{n})\overline{B_{1}}(x_{n}-x_{1})\cdots \overline{B_{n-1}}(x_{n}-x_{n-1}).
\end{align*}
\end{thm}
\noindent{\it Proof.} It is clear that the result holds for $x_{n}=0.$ We now consider the case $x_{n}\neq 0.$
By \eqref{pf31-1},
\begin{equation}\label{pf31-2}
  \sum_{\chi_{n}}{A\chi_{1}\cdots\chi_{n}\choose C\chi_{1}\cdots\chi_{n}}\chi_{n}(x_{n})=(q-1)\overline{C}\overline{\chi_{1}}\cdots\overline{\chi_{n-1}}(x_{n})\overline{A}C(1-x_{n}),
\end{equation}
which, by \eqref{pf31-1}, implies that
\begin{align*}
  &\sum_{\chi_{1},\cdots, \chi_{n-1}}{B_{1}\chi_{1}\choose \chi_{1}}\cdots {B_{n-1}\chi_{n-1}\choose \chi_{n-1}}\chi_{1}(x_{1})\cdots \chi_{n-1}(x_{n-1})\sum_{\chi_{n}}{A\chi_{1}\cdots\chi_{n}\choose C\chi_{1}\cdots\chi_{n}}\chi_{n}(x_{n})\\
  &=(q-1)\overline{C}(x_{n})\overline{A}C(1-x_{n})\sum_{\chi_{1}}{B_{1}\chi_{1}\choose \chi_{1}}\chi_{1}\left(\frac{x_{1}}{x_{n}}\right)\cdots \sum_{\chi_{n-1}}{B_{n-1}\chi_{n-1}\choose \chi_{n-1}}\chi_{n-1}\left(\frac{x_{n-1}}{x_{n}}\right)\\
  &=(q-1)^n \varepsilon (x_{1}\cdots x_{n-1})B_{1}\cdots B_{n-1}\overline{C}(x_{n})\overline{A}C(1-x_{n})\overline{B_{1}}(x_{n}-x_{1})\cdots \overline{B_{n-1}}(x_{n}-x_{n-1}).
\end{align*}
This, together with Theorem \ref{t1} and \eqref{f4}, gives
\begin{align*}
  &F^{(n)}_{D}\left(\begin{matrix}
A; B_{1},\cdots,B_{n-1},\varepsilon \\
 C \end{matrix}
\bigg| x_{1},\cdots,x_{n}\right)\\
&=\frac{1}{(q-1)^n}\sum_{\chi_{1},\cdots, \chi_{n}}{A\chi_{1}\cdots\chi_{n}\choose C\chi_{1}\cdots\chi_{n}}{B_{1}\chi_{1}\choose \chi_{1}}\cdots {B_{n-1}\chi_{n-1}\choose \chi_{n-1}}{\chi_{n}\choose \chi_{n}}\chi_{1}(x_{1})\cdots \chi_{n}(x_{n})\\
&=-\frac{1}{(q-1)^n}\sum_{\chi_{1},\cdots, \chi_{n-1}}{B_{1}\chi_{1}\choose \chi_{1}}\cdots {B_{n-1}\chi_{n-1}\choose \chi_{n-1}}\chi_{1}(x_{1})\cdots \chi_{n-1}(x_{n-1})\sum_{\chi_{n}}{A\chi_{1}\cdots\chi_{n}\choose C\chi_{1}\cdots\chi_{n}}\chi_{n}(x_{n})\\
&~+\frac{1}{(q-1)^{n-1}}\sum_{\chi_{1},\cdots, \chi_{n-1}}{A\chi_{1}\cdots\chi_{n-1}\choose C\chi_{1}\cdots\chi_{n-1}}{B_{1}\chi_{1}\choose \chi_{1}}\cdots {B_{n-1}\chi_{n-1}\choose \chi_{n-1}}\chi_{1}(x_{1})\cdots \chi_{n-1}(x_{n-1})\\
&=F^{(n-1)}_{D}\left(\begin{matrix}
A; B_{1},\cdots,B_{n-1} \\
 C \end{matrix}
\bigg| x_{1},\cdots,x_{n-1}\right)\\
&~-\varepsilon (x_{1}\cdots x_{n-1})B_{1}\cdots B_{n-1}\overline{C}(x_{n})\overline{A}C(1-x_{n})\overline{B_{1}}(x_{n}-x_{1})\cdots \overline{B_{n-1}}(x_{n}-x_{n-1}).
\end{align*}
This finishes the proof of Theorem \ref{t31}. \qed

When $n=2,$  Theorem \ref{t31} reduces to \cite [Theorem 3.1]{LLM}.  When $n=1,$ Theorem \ref{t31} reduces to \cite [Corollary 3.16, (i)]{Gr}.

It is easily seen from the definition of the Lauricella hypergeometric series $F^{(n)}_{D}$ that
\begin{equation*}
F^{(n)}_{D}\left(\begin{matrix}
a; b_{1},\cdots,b_{n} \\
 a \end{matrix}
\bigg| x_{1},\cdots,x_{n}  \right)=(1-x_{1})^{-b_{1}}\cdots(1-x_{n})^{-b_{n}}.
\end{equation*}
We also deduce the finite field analogue for the above formula.
\begin{thm}\label{t32}For $A,B_{1},\cdots,B_{n}\in \widehat{\mathbb{F}}_{q}$ and $x_{1},\cdots,x_{n}\in \mathbb{F}_{q},$ we have
\begin{align*}
  F^{(n)}_{D}\left(\begin{matrix}
A; B_{1},\cdots,B_{n} \\
 A \end{matrix}
\bigg| x_{1},\cdots,x_{n}\right)&=-\varepsilon(x_{1}\cdots x_{n})\overline{B_{1}}(1-x_{1})\cdots \overline{B_{n}}(1-x_{n})\\
&~+B_{n}(-1)\overline{A}(x_{n})F^{(n-1)}_{D}\left(\begin{matrix}
A; B_{1},\cdots,B_{n-1} \\
 A\overline{B_{n}} \end{matrix}
\bigg| \frac{x_{1}}{x_{n}},\cdots,\frac{x_{n-1}}{x_{n}}\right).
\end{align*}
\end{thm}
\noindent{\it Proof.} It follows from Theorem \ref{t1}, \eqref{f4}, \eqref{f3} and \eqref{pf31-1} that
\begin{align*}
  &F^{(n)}_{D}\left(\begin{matrix}
A; B_{1},\cdots,B_{n} \\
 A \end{matrix}
\bigg| x_{1},\cdots,x_{n}\right)\\
&=\frac{1}{(q-1)^n}\sum_{\chi_{1},\cdots, \chi_{n}}{A\chi_{1}\cdots\chi_{n}\choose A\chi_{1}\cdots\chi_{n}}{B_{1}\chi_{1}\choose \chi_{1}}\cdots {B_{n}\chi_{n}\choose \chi_{n}}\chi_{1}(x_{1})\cdots \chi_{n}(x_{n})\\
&=-\frac{1}{(q-1)^n}\sum_{\chi_{1}}{B_{1}\chi_{1}\choose \chi_{1}}\chi_{1}(x_{1})\cdots \sum_{\chi_{n}}{B_{n}\chi_{n}\choose \chi_{n}}\chi_{n}(x_{n})\\
&~+\frac{1}{(q-1)^{n-1}}\sum_{A\chi_{1}\cdots\chi_{n}=\varepsilon}{B_{1}\chi_{1}\choose \chi_{1}}\cdots {B_{n}\chi_{n}\choose \chi_{n}}\chi_{1}(x_{1})\cdots \chi_{n}(x_{n})\\
&=-\varepsilon(x_{1}\cdots x_{n})\overline{B_{1}}(1-x_{1})\cdots \overline{B_{n}}(1-x_{n})+\frac{B_{n}(-1)\overline{A}(x_{n})}{(q-1)^{n-1}}\\
&~\cdot\sum_{\chi_{1},\cdots, \chi_{n-1}}{A\chi_{1}\cdots\chi_{n-1}\choose A\overline{B_{n}}\chi_{1}\cdots\chi_{n-1}}{B_{1}\chi_{1}\choose \chi_{1}}\cdots {B_{n-1}\chi_{n-1}\choose \chi_{n-1}}\chi_{1}\left(\frac{x_{1}}{x_{n}}\right)\cdots \chi_{n-1}\left(\frac{x_{n-1}}{x_{n}}\right)\\
&=-\varepsilon(x_{1}\cdots x_{n})\overline{B_{1}}(1-x_{1})\cdots \overline{B_{n}}(1-x_{n})\\
&~+B_{n}(-1)\overline{A}(x_{n})F^{(n-1)}_{D}\left(\begin{matrix}
A; B_{1},\cdots,B_{n-1} \\
 A\overline{B_{n}} \end{matrix}
\bigg| \frac{x_{1}}{x_{n}},\cdots,\frac{x_{n-1}}{x_{n}}\right).
\end{align*}
This concludes the proof of Theorem \ref{t32}. \qed

Setting $n=2,$ we obtain the following result relating the Appell series $F_{1}$ over finite fields to the Gaussian hypergeometric series ${}_{2}F_{1}.$
\begin{cor}For  $A,B,B'\in \widehat{\mathbb{F}}_{q}$ and $x, y\in \mathbb{F}_{q},$ we have
\begin{align*}
F_{1}(A;B,B';A; x,y)&=-\varepsilon(xy)\overline{B}(1-x)\overline{B'}(1-y)+B'(-1)\overline{A}(y){}_{2}F_{1}\left(\begin{matrix}
B,A \\
 A\overline{B'}\end{matrix}
\bigg| \frac{x}{y}\right)\\
&=-\varepsilon(xy)\overline{B}(1-x)\overline{B'}(1-y)+B(-1)\overline{A}(x){}_{2}F_{1}\left(\begin{matrix}
B',A \\
 A\overline{B}\end{matrix}
\bigg| \frac{y}{x}\right).
\end{align*}
\end{cor}
When $n=1,$ Theorem \ref{t32} reduces to \cite [Corollary 3.16, (iv)]{Gr}.

The following theorem involves a transformation formula for the  Lauricella hypergeometric series over finite fields.
\begin{thm}\label{t3-6}For $A,B_{1},\cdots,B_{n},C\in \widehat{\mathbb{F}}_{q},~x_{1},\cdots,x_{n}\in \mathbb{F}_{q},$ we have
\begin{align*}
&\varepsilon((1-x_{1})\cdots(1-x_{n}))F^{(n)}_{D}\left(\begin{matrix}
A; B_{1},\cdots,B_{n} \\
 C \end{matrix}
\bigg| x_{1},\cdots,x_{n}\right)\\
&=\varepsilon(x_{1}\cdots x_{n})B_{1}\cdots B_{n}(-1)F^{(n)}_{D}\left(\begin{matrix}
A; B_{1},\cdots,B_{n} \\
 AB_{1}\cdots B_{n}\overline{C} \end{matrix}
\bigg| 1-x_{1},\cdots,1-x_{n}\right).
\end{align*}
\end{thm}
\noindent{\it Proof.} Making the substitution $u=\frac{v}{v-1}$ in the definition of  the  Lauricella hypergeometric series over finite fields gives
\begin{align*}
  &\varepsilon((1-x_{1})\cdots(1-x_{n}))F^{(n)}_{D}\left(\begin{matrix}
A; B_{1},\cdots,B_{n} \\
 C \end{matrix}
\bigg| x_{1},\cdots,x_{n}\right)\\
&=\varepsilon(x_{1}\cdots x_{n}(1-x_{1})\cdots(1-x_{n}))C(-1)\\
&~\cdot \sum_{v}A(v)B_{1}\cdots B_{n}\overline{C}(1-v)\overline{B_{1}}(1-(1-x_{1})v)\cdots \overline{B_{n}}(1-(1-x_{n})v)\\
&=\varepsilon(x_{1}\cdots x_{n})B_{1}\cdots B_{n}(-1)F^{(n)}_{D}\left(\begin{matrix}
A; B_{1},\cdots,B_{n} \\
 AB_{1}\cdots B_{n}\overline{C} \end{matrix}
\bigg| 1-x_{1},\cdots,1-x_{n}\right).
\end{align*}
This completes the proof of Theorem \ref{t3-6}.\qed

Taking $n=2$ in Theorem \ref{t3-6}, we deduce a transformation formula for  the Appell series $F_{1}$ over finite fields.

\begin{cor}For $A,B,B',C\in \widehat{\mathbb{F}}_{q},~x,y\in \mathbb{F}_{q},$ we have
\begin{align*}
&\varepsilon((1-x)(1-y))F_{1}(A; B,B';
 C ; x,y)=\varepsilon(xy)BB'(-1)F_{1}(
A; B,B';AB B'\overline{C};1-x,1-y).
\end{align*}
\end{cor}
When $n=1$ and $x=x_{1}\in \mathbb{F}_{q}\backslash \{0,1\},$  Theorem \ref{t3-6} reduces to \cite [Theorem 4.4, (i)]{Gr}.

From the integral representation for the Lauricella hypergeometric series $F^{(n)}_{D}$  we can easily obtain
\begin{align*}
  &F^{(n)}_{D}\left(\begin{matrix}
a; b_{1},\cdots,b_{n} \\
 c \end{matrix}
\bigg| x_{1},\cdots,x_{n}  \right)\\
&~=(1-x_{1})^{-b_{1}}\cdots (1-x_{n})^{-b_{n}} F^{(n)}_{D}\left(\begin{matrix}
c-a; b_{1},\cdots,b_{n} \\
 c \end{matrix}
\bigg| \frac{x_{1}}{x_{1}-1},\cdots,\frac{x_{n}}{x_{n}-1}  \right).
\end{align*}
We give a transformation formula for the Lauricella hypergeometric series over finite fields which can be regarded as the finite field analogue for the above identity.
\begin{thm}\label{t33}For  $A,B_{1},\cdots,B_{n},C\in \widehat{\mathbb{F}}_{q}$ and $x_{1},\cdots,x_{n}\in \mathbb{F}_{q}\backslash \{1\},$  we have
\begin{align*}
  &F^{(n)}_{D}\left(\begin{matrix}
A; B_{1},\cdots,B_{n} \\
 C\end{matrix}
\bigg| x_{1},\cdots,x_{n}\right)\\
&=C(-1)\overline{B_{1}}(1-x_{1})\cdots \overline{B_{n}}(1-x_{n})F^{(n)}_{D}\left(\begin{matrix}
\overline{A}C; B_{1},\cdots,B_{n} \\
 C\end{matrix}
\bigg| \frac{x_{1}}{x_{1}-1},\cdots,\frac{x_{n}}{x_{n}-1}\right).
\end{align*}
\end{thm}
\noindent{\it Proof.}  The result follows from the definition of the Lauricella hypergeometric series over finite fields and Making the substitution $u=1-v.$ \qed

When $n=2,$  Theorem \ref{t33} reduces to \cite [Theorem 3.2, (3.6)]{LLM}.  When $n=1,$  Theorem \ref{t33} reduces to \cite [Theorem 4.4, (ii)]{Gr} for $x\neq 1.$

\begin{thm}\label{t34}For $A,B_{1},\cdots,B_{n},C\in \widehat{\mathbb{F}}_{q},~x_{1},\cdots,x_{n-1}\in \mathbb{F}_{q}$ and $x_{n}\in \mathbb{F}_{q}\backslash \{1\},$  we have
\begin{align*}
  &\varepsilon((x_{n}-x_{1})\cdots (x_{n}-x_{n-1}))F^{(n)}_{D}\left(\begin{matrix}
A; B_{1},\cdots,B_{n} \\
 C\end{matrix}
\bigg| x_{1},\cdots,x_{n}\right)\\
&=\varepsilon(x_{1}\cdots x_{n-1})\overline{A}(1-x_{n})F^{(n)}_{D}\left(\begin{matrix}
A; B_{1},\cdots,B_{n-1}, \overline{B_{1}}\cdots \overline{B_{n}}C \\
 C\end{matrix}
\bigg| \frac{x_{n}-x_{1}}{x_{n}-1},\cdots,\frac{x_{n}-x_{n-1}}{x_{n}-1}, \frac{x_{n}}{x_{n}-1}\right).
\end{align*}
\end{thm}
\noindent{\it Proof.} Making the substitution $u=\frac{v}{1-x_{n}+x_{n}v}$ in the the definition of the Lauricella hypergeometric series over finite fields, we have
\begin{align*}
  &\varepsilon((x_{n}-x_{1})\cdots (x_{n}-x_{n-1}))F^{(n)}_{D}\left(\begin{matrix}
A; B_{1},\cdots,B_{n} \\
 C\end{matrix}
\bigg| x_{1},\cdots,x_{n}\right)\\
&=\varepsilon(x_{1}\cdots x_{n}(x_{n}-x_{1})\cdots (x_{n}-x_{n-1}))AC(-1)\sum_{u}A(u)\overline{A}C(1-u)\overline{B_{1}}(1-x_{1}u)\cdots \overline{B_{n}}(1-x_{n}u)\\
&=\varepsilon(x_{1}\cdots x_{n}(x_{n}-x_{1})\cdots (x_{n}-x_{n-1}))AC(-1)\overline{A}(1-x_{n})\\
&~\cdot\sum_{v}A(v)\overline{A}C(1-v)\overline{B_{1}}\left(1-\frac{x_{n}-x_{1}}{x_{n}-1}v\right)\cdots\overline{B_{n-1}}\left(1-\frac{x_{n}-x_{n-1}}{x_{n}-1}v\right)
B_{1}\cdots B_{n}\overline{C}\left(1-\frac{x_{n}}{x_{n}-1}v\right)\\
&=\varepsilon(x_{1}\cdots x_{n-1})\overline{A}(1-x_{n})F^{(n)}_{D}\left(\begin{matrix}
A; B_{1},\cdots,B_{n-1}, \overline{B_{1}}\cdots \overline{B_{n}}C \\
 C\end{matrix}
\bigg| \frac{x_{n}-x_{1}}{x_{n}-1},\cdots,\frac{x_{n}-x_{n-1}}{x_{n}-1}, \frac{x_{n}}{x_{n}-1}\right),
\end{align*}
from which we complete the proof of Theorem \ref{t34}.\qed

From Theorem \ref{t34} and Theorem \ref{t31}, we can easily obtain the following reduction formula for the Lauricella hypergeometric series over finite fields.
\begin{cor} \label{c35}For $A,B_{1},\cdots,B_{n},\in \widehat{\mathbb{F}}_{q},~x_{1},\cdots,x_{n-1}\in \mathbb{F}_{q}$ and $x_{n}\in \mathbb{F}_{q}\backslash \{1\},$  we have
\begin{align*}
  &\varepsilon((x_{n}-x_{1})\cdots (x_{n}-x_{n-1}))F^{(n)}_{D}\left(\begin{matrix}
A; B_{1},\cdots,B_{n} \\
 B_{1}\cdots B_{n}\end{matrix}
\bigg| x_{1},\cdots,x_{n}\right)\\
&=\varepsilon(x_{1}\cdots x_{n})\overline{A}(1-x_{n})F^{(n-1)}_{D}\left(\begin{matrix}
A; B_{1},\cdots,B_{n-1} \\
 B_{1}\cdots B_{n}\end{matrix}
\bigg| \frac{x_{n}-x_{1}}{x_{n}-1},\cdots,\frac{x_{n}-x_{n-1}}{x_{n}-1}\right)\\
&~-\varepsilon((x_{n}-x_{1})\cdots(x_{n}-x_{n-1}))\overline{B_{1}}(-x_{1})\cdots \overline{B_{n}}(-x_{n}).
\end{align*}
\end{cor}
Actually, the formula in Corollary \ref{c35} can be considered as a finite field analogue for  the following reduction formula on the Lauricella hypergeometric series (see G. Mingari Scarpello and D. Ritelli \cite{MR}):
\begin{equation*}
F^{(n)}_{D}\left(\begin{matrix}
a; b_{1},\cdots,b_{n} \\
 b_{1}+\cdots+b_{n} \end{matrix}
\bigg| x_{1},\cdots,x_{n}  \right)=\frac{1}{(1-x_{n})^{a}}F^{(n-1)}_{D}\left(\begin{matrix}
a; b_{1},\cdots,b_{n-1} \\
 b_{1}+\cdots+b_{n} \end{matrix}
\bigg| \frac{x_{1}-x_{n}}{1-x_{n}},\cdots,\frac{x_{n-1}-x_{n}}{1-x_{n}}\right).
\end{equation*}

When $n=2,$ Theorem \ref{t34} reduces to \cite [Theorem 3.2, (3.7) and (3.9)]{LLM}. When $n=1,$ Theorem \ref{t34} reduces to \cite [Theorem 4.4, (iii)]{Gr} for $x\neq 1.$
\begin{thm}\label{t36} For $A,B_{1},\cdots,B_{n},C\in \widehat{\mathbb{F}}_{q}$ and $x_{1},\cdots,x_{n}\in \mathbb{F}_{q}\backslash \{1\},$  we have
\begin{align*}
  &\varepsilon((x_{n}-x_{1})\cdots (x_{n}-x_{n-1}))F^{(n)}_{D}\left(\begin{matrix}
A; B_{1},\cdots,B_{n} \\
 C\end{matrix}
\bigg| x_{1},\cdots,x_{n}\right)\\
&=\varepsilon(x_{1}\cdots x_{n-1})C(-1)\overline{A}\overline{B_{n}}C(1-x_{n})\overline{B_{1}}(1-x_{1})\cdots \overline{B_{n-1}}(1-x_{n-1})\\
&~\cdot F^{(n)}_{D}\left(\begin{matrix}
\overline{A}C; B_{1},\cdots,B_{n-1},  \overline{B_{1}}\cdots \overline{B_{n}}C\\
 C\end{matrix}
\bigg| \frac{x_{n}-x_{1}}{1-x_{1}},\cdots,\frac{x_{n}-x_{n-1}}{1-x_{n-1}}, x_{n}\right)
\end{align*}
\end{thm}
\noindent{\it Proof.} Making another substitution $u=\frac{1-v}{1-vx_{n}}$ in the the definition of the Lauricella hypergeometric series over finite fields,  we get
\begin{align*}
  &\varepsilon((x_{n}-x_{1})\cdots (x_{n}-x_{n-1}))F^{(n)}_{D}\left(\begin{matrix}
A; B_{1},\cdots,B_{n} \\
 C\end{matrix}
\bigg| x_{1},\cdots,x_{n}\right)\\
&=\varepsilon(x_{1}\cdots x_{n}(x_{n}-x_{1})\cdots (x_{n}-x_{n-1}))AC(-1)\overline{A}\overline{B_{n}}C(1-x_{n})\overline{B_{1}}(1-x_{1})\cdots \overline{B_{n-1}}(1-x_{n-1})\\
&~\cdot \sum_{v}\overline{A}C(v)A(1-v)\overline{B_{1}}\left(1-\frac{x_{n}-x_{1}}{1-x_{1}}v\right)\cdots \overline{B_{n-1}}\left(1-\frac{x_{n}-x_{n-1}}{1-x_{n-1}}v\right)B_{1}\cdots B_{n}\overline{C}(1-x_{n}v)\\
&=\varepsilon(x_{1}\cdots x_{n-1})C(-1)\overline{A}\overline{B_{n}}C(1-x_{n})\overline{B_{1}}(1-x_{1})\cdots \overline{B_{n-1}}(1-x_{n-1})\\
&~\cdot F^{(n)}_{D}\left(\begin{matrix}
\overline{A}C; B_{1},\cdots,B_{n-1},  \overline{B_{1}}\cdots \overline{B_{n}}C\\
 C\end{matrix}
\bigg| \frac{x_{n}-x_{1}}{1-x_{1}},\cdots,\frac{x_{n}-x_{n-1}}{1-x_{n-1}}, x_{n}\right).
\end{align*}
This completes the proof of Theorem \ref{t36}.\qed

Similarly, we can get another reduction formula.
\begin{cor}\label{c37} For $A,B_{1},\cdots,B_{n}\in \widehat{\mathbb{F}}_{q}$ and $x_{1},\cdots,x_{n}\in \mathbb{F}_{q}\backslash \{1\},$  we have
\begin{align*}
  &\varepsilon((x_{n}-x_{1})\cdots (x_{n}-x_{n-1}))F^{(n)}_{D}\left(\begin{matrix}
A; B_{1},\cdots,B_{n} \\
 B_{1}\cdots B_{n}\end{matrix}
\bigg| x_{1},\cdots,x_{n}\right)\\
&=\varepsilon(x_{1}\cdots x_{n})B_{1}\cdots B_{n}(-1)\overline{A}B_{1}\cdots B_{n-1}(1-x_{n})\overline{B_{1}}(1-x_{1})\cdots \overline{B_{n-1}}(1-x_{n-1})\\
&~\cdot F^{(n-1)}_{D}\left(\begin{matrix}
\overline{A}B_{1}\cdots B_{n}; B_{1},\cdots,B_{n-1}\\
 B_{1}\cdots B_{n}\end{matrix}
\bigg| \frac{x_{n}-x_{1}}{1-x_{1}},\cdots,\frac{x_{n}-x_{n-1}}{1-x_{n-1}}\right)\\
&~-\varepsilon((x_{n}-1)(x_{n}-x_{1})\cdots (x_{n}-x_{n-1}))\overline{B_{1}}(-x_{1})\cdots\overline{B_{n}}(-x_{n}).
\end{align*}
\end{cor}
When $n=2,$ Theorem \ref{t36} reduces to \cite [(3.8) and (3.10)]{LLM}. When $n=1,$  Theorem \ref{t36} reduces to \cite [Theorem 4.4, (iv)]{Gr} for $x\neq 1.$

We also  give some evaluations for the Lauricella hypergeometric series over finite fields.

From the definition of the Lauricella hypergeometric series over finite fields,  we can easily deduce the following results.
\begin{thm}\label{t6-1}For $A,B_{1},\cdots,B_{n},C\in \widehat{\mathbb{F}}_{q}$ and $x, x_{1},\cdots,x_{n-1}\in \mathbb{F}_{q},$ we have
\begin{align*}
F^{(n)}_{D}\left(\begin{matrix}
A; B_{1},\cdots,B_{n} \\
 C \end{matrix}
\bigg| x,\cdots,x\right)&={}_{2}F_1 \left(\begin{matrix}
 B_{1}\cdots B_{n},A \\
C \end{matrix}
\bigg| x \right),\\
  F^{(n)}_{D}\left(\begin{matrix}
A; B_{1},\cdots,B_{n} \\
 C \end{matrix}
\bigg| x_{1},\cdots,x_{n-1}, 1\right)&=B_{n}(-1)F^{(n-1)}_{D}\left(\begin{matrix}
A; B_{1},\cdots,B_{n-1} \\
 \overline{B_{n}}C \end{matrix}
\bigg| x_{1},\cdots,x_{n-1}\right).
\end{align*}
In particular,
\begin{align}\label{t61-1}
  F^{(n)}_{D}\left(\begin{matrix}
A; B_{1},\cdots,B_{n} \\
 C \end{matrix}
\bigg| 1,\cdots, 1\right)=B_{1}\cdots B_{n}(-1){A\choose \overline{B_{1}}\cdots \overline{B_{n}}C}.
\end{align}
\end{thm}
When $n=1,$ \eqref{t61-1} reduces to \cite [Theorem 4.9]{Gr}:
\begin{equation}\label{t61-2}
{}_{2}F_{1}\left(\begin{matrix}
A, B \\
 C \end{matrix}
\bigg| 1\right)=A(-1){B\choose \overline{A}C}.
\end{equation}

\begin{cor}\label{c6-2} For $A,B_{1},\cdots,B_{n}\in \widehat{\mathbb{F}}_{q}$ and $x\in \mathbb{F}_{q},$ we have
\begin{align*}
  F^{(n)}_{D}\left(\begin{matrix}
A; B_{1},\cdots,B_{n} \\
 A \end{matrix}
\bigg| x,\cdots,x\right)=-\varepsilon(x)\overline{B_{1}}\cdots \overline{B_{n}}(1-x)+B_{1}\cdots B_{n}(-1)\overline{A}(x){A\choose B_{1}\cdots B_{n}}.
\end{align*}
\end{cor}
\noindent{\it Proof.}  We can take $x_{1}=\cdots=x_{n}=x$ in Theorem \ref{t32} and using \eqref{t61-1} to get the result. Alternatively, we take $C=A$ in the first identity of Theorem \ref{t6-1} and use \eqref{f2}--\eqref{f4} and \eqref{pf31-1}:
\begin{align*}
  F^{(n)}_{D}\left(\begin{matrix}
A; B_{1},\cdots,B_{n} \\
 A \end{matrix}
\bigg| x,\cdots,x\right)&={}_{2}F_1 \left(\begin{matrix}
 B_{1}\cdots B_{n},A \\
A \end{matrix}
\bigg| x \right)\\
&=\frac{1}{q-1}\sum_{\chi}{B_{1}\cdots B_{n}\chi \choose \chi}\chi(x)+\overline{A}(x){\overline{A}B_{1}\cdots B_{n}\choose \overline{A}}\\
&=-\varepsilon(x)\overline{B_{1}}\cdots \overline{B_{n}}(1-x)+B_{1}\cdots B_{n}(-1)\overline{A}(x){A\choose B_{1}\cdots B_{n}}
\end{align*}
to complete the proof of Corollary \ref{c6-2}. \qed

Setting $n=2$ in Corollary \ref{c6-2}, we are led to
\begin{cor}For $A,B,B'\in \widehat{\mathbb{F}}_{q}$ and $x\in \mathbb{F}_{q},$ we have
\begin{align*}
  F_{1}(A; B,B';A; x,x)=-\varepsilon(x)\overline{B}\overline{B'}(1-x)+BB'(-1)\overline{A}(x){A\choose B B'}.
\end{align*}
\end{cor}
From Theorem \ref{t6-1} we can obtain another result.
\begin{cor}\label{c6-3}For $A,B_{1},\cdots,B_{n}\in \widehat{\mathbb{F}}_{q}$ and $x\in \mathbb{F}_{q},$ we have
\begin{align*}
  F^{(n)}_{D}\left(\begin{matrix}
A; B_{1},\cdots,B_{n} \\
 B_{1}\cdots B_{n} \end{matrix}
\bigg| x,\cdots,x\right)&={A\choose B_{1}\cdots B_{n}}\varepsilon(x)\overline{A}(1-x)-\overline{B_{1}}\cdots \overline{B_{n}}(-x)\\
&~+(q-1)B_{1}\cdots B_{n}(-1)\delta(1-x)\delta(A).
\end{align*}
\end{cor}
\noindent{\it Proof.} It follows from \cite [Corollary 3.16, (iii)]{Gr} that
\begin{align*}
{}_{2}F_1 \left(\begin{matrix}
A, B \\
A \end{matrix}
\bigg| x \right)={B\choose A}\varepsilon(x)\overline{B}(1-x)-\overline{A}(-x)+(q-1)A(-1)\delta(1-x)\delta(B)
\end{align*}
for $A,B\in \widehat{\mathbb{F}}_{q}$ and $x\in \mathbb{F}_{q}.$ Then, putting $C=B_{1}\cdots B_{n}$ in the first identity of Theorem \ref{t6-1},  we have
\begin{align*}
  F^{(n)}_{D}\left(\begin{matrix}
A; B_{1},\cdots,B_{n} \\
 B_{1}\cdots B_{n} \end{matrix}
\bigg| x,\cdots,x\right)&={}_{2}F_1 \left(\begin{matrix}
 B_{1}\cdots B_{n},A \\
B_{1}\cdots B_{n} \end{matrix}
\bigg| x \right)\\
&={A\choose B_{1}\cdots B_{n}}\varepsilon(x)\overline{A}(1-x)-\overline{B_{1}}\cdots \overline{B_{n}}(-x)\\
&~+(q-1)B_{1}\cdots B_{n}(-1)\delta(1-x)\delta(A),
\end{align*}
from which the result follows.\qed

Letting $n=2$ in Corollary \ref{c6-3} gives
\begin{cor}For $A,B,B'\in \widehat{\mathbb{F}}_{q}$ and $x\in \mathbb{F}_{q},$ we have
\begin{align*}
  F_{1}(A; B,B';B B';x,x)={A\choose BB'}\varepsilon(x)\overline{A}(1-x)-\overline{B}\overline{B'}(-x)
+(q-1)B B'(-1)\delta(1-x)\delta(A).
\end{align*}
\end{cor}
\section{Generating functions}
In this section, we establish several  generating functions for the Lauricella hypergeometric series over finite fields.

\begin{thm}\label{t4-1}For $A,B_{1},\cdots,B_{n},C\in \widehat{\mathbb{F}}_{q}, x_{1},\cdots,x_{n}\in \mathbb{F}_{q}$ and $t\in \mathbb{F}_{q}\backslash \{1\},$ we have
\begin{align*}
  &\sum_{\theta} {A\overline{C}\theta\choose\theta}F^{(n)}_{D}\left(\begin{matrix}
A\theta; B_{1},\cdots,B_{n} \\
 C \end{matrix}
\bigg| x_{1},\cdots,x_{n}\right)\theta(t)\\
&=\varepsilon(t)\overline{A}(1-t)F^{(n)}_{D}\left(\begin{matrix}
A; B_{1},\cdots,B_{n} \\
 C \end{matrix}
\bigg| \frac{x_{1}}{1-t},\cdots,\frac{x_{n}}{1-t}\right)\\
&~-\varepsilon(x_{1}\cdots x_{n})\overline{A}C(-t)\overline{B_{1}}(1-x_{1})\cdots \overline{B_{n}}(1-x_{n}).
\end{align*}
\end{thm}
\noindent{\it Proof.} Making the substitution $u=\frac{v}{1-t},$ we have
\begin{align*}
  &\varepsilon(t x_{1}\cdots x_{n})AC(-1)\sum_{u\neq 1}A(u)\overline{A}C(1-u+ut)\overline{B_{1}}(1-x_{1}u)\cdots \overline{B_{n}}(1-x_{n}u)\\
  &=\varepsilon(t x_{1}\cdots x_{n})AC(-1)\sum_{u}A(u)\overline{A}C(1-u+ut)\overline{B_{1}}(1-x_{1}u)\cdots \overline{B_{n}}(1-x_{n}u)\\
  &~-\varepsilon(x_{1}\cdots x_{n})\overline{A}C(-t)\overline{B_{1}}(1-x_{1})\cdots \overline{B_{n}}(1-x_{n})\\
  &=\varepsilon(t x_{1}\cdots x_{n})AC(-1)\overline{A}(1-t)\sum_{v}A(v)\overline{A}C(1-v)\overline{B_{1}}\left(1-\frac{x_{1}}{1-t}v\right)\cdots \overline{B_{n}}\left(1-\frac{x_{n}}{1-t}v\right)\\
  &~-\varepsilon(x_{1}\cdots x_{n})\overline{A}C(-t)\overline{B_{1}}(1-x_{1})\cdots \overline{B_{n}}(1-x_{n})\\
  &=\varepsilon(t)\overline{A}(1-t)F^{(n)}_{D}\left(\begin{matrix}
A; B_{1},\cdots,B_{n} \\
 C \end{matrix}
\bigg| \frac{x_{1}}{1-t},\cdots,\frac{x_{n}}{1-t}\right)\\
&~-\varepsilon(x_{1}\cdots x_{n})\overline{A}C(-t)\overline{B_{1}}(1-x_{1})\cdots \overline{B_{n}}(1-x_{n}).
\end{align*}
This combines the binomial theorem over finite fields  to yield
\begin{align*}
&\sum_{\theta} {A\overline{C}\theta\choose\theta}F^{(n)}_{D}\left(\begin{matrix}
A\theta; B_{1},\cdots,B_{n} \\
 C \end{matrix}
\bigg| x_{1},\cdots,x_{n}\right)\theta(t)\\
&=\varepsilon(x_{1}\cdots x_{n})AC(-1)\sum_{\theta,u}{A\overline{C}\theta\choose\theta}A(u)\overline{A}C(1-u)\theta(-ut)\overline{\theta}(1-u)\overline{B_{1}}(1-x_{1}u)\cdots \overline{B_{n}}(1-x_{n}u)\\
&=\varepsilon(x_{1}\cdots x_{n})AC(-1)\sum_{u\neq 1}A(u)\overline{A}C(1-u)\overline{B_{1}}(1-x_{1}u)\cdots \overline{B_{n}}(1-x_{n}u)\sum_{\theta}{A\overline{C}\theta\choose\theta}\theta\left(-\frac{ut}{1-u}\right)\\
&=\varepsilon(t x_{1}\cdots x_{n})AC(-1)\sum_{u\neq 1}A(u)\overline{A}C(1-u)\overline{A}C\left(1+\frac{ut}{1-u}\right)\overline{B_{1}}(1-x_{1}u)\cdots \overline{B_{n}}(1-x_{n}u)
\end{align*}
\begin{align*}
&=\varepsilon(t x_{1}\cdots x_{n})AC(-1)\sum_{u\neq 1}A(u)\overline{A}C(1-u+ut)\overline{B_{1}}(1-x_{1}u)\cdots \overline{B_{n}}(1-x_{n}u)\\
&=\varepsilon(t)\overline{A}(1-t)F^{(n)}_{D}\left(\begin{matrix}
A; B_{1},\cdots,B_{n} \\
 C \end{matrix}
\bigg| \frac{x_{1}}{1-t},\cdots,\frac{x_{n}}{1-t}\right)\\
&~-\varepsilon(x_{1}\cdots x_{n})\overline{A}C(-t)\overline{B_{1}}(1-x_{1})\cdots \overline{B_{n}}(1-x_{n}),
\end{align*}
which ends the proof of Theorem \ref{t4-1}.\qed

Theorem \ref{t4-1} reduces to \cite [Theorem 4.1]{LLM} when $n=2.$

Setting $n=1$ in Theorem \ref{t4-1}, we get a generating function for the Gaussian  hypergeometric series  ${}_{2}F_{1}.$
\begin{cor}For $A,B,C\in \widehat{\mathbb{F}}_{q}, x\in \mathbb{F}_{q}$ and $t\in \mathbb{F}_{q}\backslash \{1\},$ we have
\begin{align*}
  &\sum_{\theta} {B\overline{C}\theta\choose\theta}{}_{2}F_{1}\left(\begin{matrix}
A, B\theta \\
 C \end{matrix}
\bigg| x\right)\theta(t)=\varepsilon(t)\overline{B}(1-t){}_{2}F_{1}\left(\begin{matrix}
A, B\\
 C \end{matrix}
\bigg| \frac{x}{1-t}\right)-\varepsilon(x)\overline{B}C(-t)\overline{A}(1-x).
\end{align*}
\end{cor}

We also give another generating function for the Lauricella hypergeometric series over finite fields.
\begin{thm}\label{t4-2}For $A,B_{1},\cdots,B_{n},C\in \widehat{\mathbb{F}}_{q}$ and $x_{1},\cdots,x_{n}, t\in \mathbb{F}_{q},$ we have
\begin{align*}
  &\sum_{\theta}{B_{n}\theta\choose\theta}F^{(n)}_{D}\left(\begin{matrix}
A; B_{1},\cdots,B_{n-1}, B_{n}\theta \\
 C \end{matrix}
\bigg| x_{1},\cdots,x_{n}\right)\theta(t)\\
&=(q-1)\varepsilon(t)\overline{B_{n}}(1-t)F^{(n)}_{D}\left(\begin{matrix}
A; B_{1},\cdots,B_{n-1}, B_{n} \\
 C \end{matrix}
\bigg| x_{1},\cdots,x_{n-1}, \frac{x_{n}}{1-t}\right)\\
&~-(q-1)\varepsilon(x_{1}\cdots x_{n-1})\overline{B_{n}}(-t)B_{1}\cdots B_{n-1}\overline{C}(x_{n})\overline{A}C(1-x_{n})\overline{B_{1}}\left(x_{n}-x_{1}\right)\cdots \overline{B_{n-1}}\left(x_{n}-x_{n-1}\right).
\end{align*}
\end{thm}
\noindent{\it Proof.} It is obvious that the result holds for $x_{n}=0.$ We now consider the case $x_{n}\neq 0.$ It follows from \cite [Corollary 3.16, (iii)]{Gr} that
\begin{equation*}
\sum_{\theta}{B_{n}\theta\choose\theta}{B_{n}\chi_{n}\theta\choose B_{n}\theta}\theta(t)=(q-1)\left(\varepsilon(t)\overline{B_{n}}\overline{\chi_{n}}(1-t){B_{n}\chi_{n}\choose \chi_{n}}-\overline{B_{n}}(-t)\right).
\end{equation*}
Then, by \eqref{pf31-2},
\begin{align*}
&\sum_{\theta}{B_{n}\theta\choose\theta}F^{(n)}_{D}\left(\begin{matrix}
A; B_{1},\cdots,B_{n-1}, B_{n}\theta \\
 C \end{matrix}
\bigg| x_{1},\cdots,x_{n}\right)\theta(t)\\
&=\frac{1}{(q-1)^n}\sum_{\chi_{1},\cdots, \chi_{n}}{A\chi_{1}\cdots\chi_{n}\choose C\chi_{1}\cdots\chi_{n}}{B_{1}\chi_{1}\choose \chi_{1}}\cdots {B_{n-1}\chi_{n-1}\choose \chi_{n-1}}\chi_{1}(x_{1})\cdots \chi_{n}(x_{n})\\
&~\cdot  \sum_{\theta}{B_{n}\theta\choose\theta}{B_{n}\chi_{n}\theta\choose B_{n}\theta}\theta(t)
\end{align*}
\begin{align*}
&=\frac{\varepsilon(t)\overline{B_{n}}(1-t)}{(q-1)^{n-1}}\sum_{\chi_{1},\cdots, \chi_{n}}{A\chi_{1}\cdots\chi_{n}\choose C\chi_{1}\cdots\chi_{n}}{B_{1}\chi_{1}\choose \chi_{1}}\cdots {B_{n}\chi_{n}\choose \chi_{n}}\chi_{1}(x_{1})\cdots \chi_{n-1}(x_{n-1}) \chi_{n}\left(\frac{x_{n}}{1-t}\right)\\
&~-\frac{\overline{B_{n}}(-t)}{(q-1)^{n-1}}\sum_{\chi_{1},\cdots, \chi_{n-1}}{B_{1}\chi_{1}\choose \chi_{1}}\cdots {B_{n-1}\chi_{n-1}\choose \chi_{n-1}}\chi_{1}(x_{1})\cdots \chi_{n-1}(x_{n-1})\sum_{\chi_{n}}{A\chi_{1}\cdots\chi_{n}\choose C\chi_{1}\cdots\chi_{n}}\chi_{n}(x_{n})\\
&=(q-1)\varepsilon(t)\overline{B_{n}}(1-t)F^{(n)}_{D}\left(\begin{matrix}
A; B_{1},\cdots,B_{n-1}, B_{n} \\
 C \end{matrix}
\bigg| x_{1},\cdots,x_{n-1}, \frac{x_{n}}{1-t}\right)\\
&~-\frac{\overline{B_{n}}(-t)\overline{C}(x_{n})\overline{A}C(1-x_{n})}{(q-1)^{n-2}}\sum_{\chi_{1}}{B_{1}\chi_{1}\choose \chi_{1}}\chi_{1}\left(\frac{x_{1}}{x_{n}}\right)\cdots \sum_{\chi_{n-1}}{B_{n-1}\chi_{n-1}\choose \chi_{n-1}}\chi_{n-1}\left(\frac{x_{n-1}}{x_{n}}\right)\\
&=(q-1)\varepsilon(t)\overline{B_{n}}(1-t)F^{(n)}_{D}\left(\begin{matrix}
A; B_{1},\cdots,B_{n-1}, B_{n} \\
 C \end{matrix}
\bigg| x_{1},\cdots,x_{n-1}, \frac{x_{n}}{1-t}\right)\\
&~-(q-1)\varepsilon(x_{1}\cdots x_{n-1})\overline{B_{n}}(-t)\overline{C}(x_{n})\overline{A}C(1-x_{n})\overline{B_{1}}\left(1-\frac{x_{1}}{x_{n}}\right)\cdots \overline{B_{n-1}}\left(1-\frac{x_{n-1}}{x_{n}}\right),
\end{align*}
from which the result follows. This finishes the proof of Theorem \ref{t4-2}. \qed

When $n=2,$ Theorem \ref{t4-2} reduces to \cite [Theorem 4.2]{LLM}.

Taking $n=1$ in Theorem \ref{t4-2}, we can easily obtain another generating function for the Gaussian  hypergeometric series  ${}_{2}F_{1}.$

\begin{cor}For $A,B,C\in \widehat{\mathbb{F}}_{q}$ and $x,t\in \mathbb{F}_{q},$ we have
\begin{align*}
  \sum_{\theta}{A\theta\choose\theta}{}_{2}F_{1}\left(\begin{matrix}
A\theta,  B \\
 C \end{matrix}
\bigg| x\right)\theta(t)&=(q-1)\varepsilon(t)\overline{A}(1-t){}_{2}F_{1}\left(\begin{matrix}
A, B \\
 C \end{matrix}
\bigg| \frac{x}{1-t}\right)\\
&~-(q-1)\overline{A}(-t)\overline{C}(x)\overline{B}C(1-x).
\end{align*}
\end{cor}

The following theorem involves another generating function for the Lauricella hypergeometric series over finite fields.
\begin{thm}\label{t4-3}For $A,B_{1},\cdots,B_{n},C\in \widehat{\mathbb{F}}_{q}$ and  $ x_{1},\cdots,x_{n}, t\in \mathbb{F}_{q},$ we have
\begin{align*}
  &\sum_{\theta}{A\overline{C}\theta\choose\theta}F^{(n)}_{D}\left(\begin{matrix}
A; B_{1},\cdots, B_{n} \\
 C\overline{\theta} \end{matrix}
\bigg| x_{1},\cdots,x_{n}\right)\theta(t)\\
&=(q-1)\varepsilon(t)\overline{C}(1+t)F^{(n)}_{D}\left(\begin{matrix}
A; B_{1},\cdots, B_{n} \\
 C \end{matrix}
\bigg| \frac{x_{1}}{1+t},\cdots,\frac{x_{n}}{1+t}\right)\\
&~-(q-1)\overline{A}C(-t)\varepsilon(x_{1}\cdots x_{n})\overline{B_{1}}(1-x_{1})\cdots \overline{B_{n}}(1-x_{n}).
\end{align*}
\end{thm}
\noindent{\it Proof.} It is easily seen from  \cite [Corollary 3.16, (iii)]{Gr} and \eqref{f2}, \eqref{f3} that
\begin{align*}
  \sum_{\theta}{A\overline{C}\theta\choose\theta}{\overline{C}\overline{\chi_{1}}\cdots \overline{\chi_{n}}\theta \choose A\overline{C}\theta}\theta(-t)
  &=(q-1){A\chi_{1}\cdots \chi_{n}\choose C\chi_{1}\cdots \chi_{n}}AC(-1)\varepsilon(t)\overline{C}\overline{\chi_{1}}\cdots \overline{\chi_{n}}(1+t)\\
  &~-(q-1)\overline{A}C(t).
\end{align*}
Combining \eqref{f2}, \eqref{f3} and  the above identity, we obtain
\begin{align*}
  &\sum_{\theta}{A\overline{C}\theta\choose\theta}F^{(n)}_{D}\left(\begin{matrix}
A; B_{1},\cdots, B_{n} \\
 C\overline{\theta} \end{matrix}
\bigg| x_{1},\cdots,x_{n}\right)\theta(t)\\
&=\frac{AC(-1)}{(q-1)^n}\sum_{\theta,\chi_{1},\cdots, \chi_{n}}{A\overline{C}\theta\choose\theta}{\overline{C}\overline{\chi_{1}}\cdots \overline{\chi_{n}}\theta \choose A\overline{C}\theta}{B_{1}\chi_{1}\choose \chi_{1}}\cdots {B_{n}\chi_{n}\choose \chi_{n}}\chi_{1}(x_{1})\cdots \chi_{n}(x_{n})\theta(-t)\\
&=\frac{AC(-1)}{(q-1)^n}\sum_{\chi_{1},\cdots, \chi_{n}}{B_{1}\chi_{1}\choose \chi_{1}}\cdots {B_{n}\chi_{n}\choose \chi_{n}}\chi_{1}(x_{1})\cdots \chi_{n}(x_{n})\sum_{\theta}{A\overline{C}\theta\choose\theta}{\overline{C}\overline{\chi_{1}}\cdots \overline{\chi_{n}}\theta \choose A\overline{C}\theta}\theta(-t)\\
&=\frac{\varepsilon(t)\overline{C}(1+t)}{(q-1)^{n-1}}\sum_{\chi_{1},\cdots, \chi_{n}}{A\chi_{1}\cdots \chi_{n}\choose C\chi_{1}\cdots \chi_{n}}{B_{1}\chi_{1}\choose \chi_{1}}\cdots {B_{n}\chi_{n}\choose \chi_{n}}\chi_{1}\left(\frac{x_{1}}{1+t}\right)\cdots \chi_{n}\left(\frac{x_{n}}{1+t}\right)\\
&~-\frac{\overline{A}C(-t)}{(q-1)^{n-1}}\sum_{\chi_{1}}{B_{1}\chi_{1}\choose \chi_{1}}\chi_{1}(x_{1})\cdots \sum_{\chi_{n}}{B_{n}\chi_{n}\choose \chi_{n}} \chi_{n}(x_{n})\\
&=(q-1)\varepsilon(t)\overline{C}(1+t)F^{(n)}_{D}\left(\begin{matrix}
A; B_{1},\cdots, B_{n} \\
 C \end{matrix}
\bigg| \frac{x_{1}}{1+t},\cdots,\frac{x_{n}}{1+t}\right)\\
&~-(q-1)\overline{A}C(-t)\varepsilon(x_{1}\cdots x_{n})\overline{B_{1}}(1-x_{1})\cdots \overline{B_{n}}(1-x_{n}).
\end{align*}
This concludes the proof of Theorem \ref{t4-3}.\qed

Putting $n=2$ in Theorem \ref{t4-3}, we get the following result which is a generating function for the finite field analogue of the Appell series $F_{1}.$
\begin{cor}For $A,B,B',C\in \widehat{\mathbb{F}}_{q}$ and  $ x,y,t\in \mathbb{F}_{q},$ we have
\begin{align*}
  \sum_{\theta}{A\overline{C}\theta\choose\theta}F_{1}(
A; B, B' ;C\overline{\theta} ; x,y)\theta(t)&=(q-1)\varepsilon(t)\overline{C}(1+t)F_{1}\left(A; B_, B' ;  C ; \frac{x}{1+t},\frac{y}{1+t}\right)\\
&~-(q-1)\overline{A}C(-t)\varepsilon(xy)\overline{B}(1-x)\overline{B'}(1-y).
\end{align*}
\end{cor}
Letting  $n=1$ in Theorem \ref{t4-3} yields a generating function for the Gaussian  hypergeometric series  ${}_{2}F_{1}.$
\begin{cor}For $A,B,C\in \widehat{\mathbb{F}}_{q}$ and  $x,t\in \mathbb{F}_{q},$ we have
\begin{align*}
\sum_{\theta}{B\overline{C}\theta\choose\theta}{}_{2}F_{1}\left(\begin{matrix}
A, B \\
 C\overline{\theta} \end{matrix}
\bigg| x\right)\theta(t)
&=(q-1)\varepsilon(t)\overline{C}(1+t){}_{2}F_{1}\left(\begin{matrix}
A, B \\
 C \end{matrix}
\bigg| \frac{x}{1+t}\right)\\
&~-(q-1)\overline{B}C(-t)\varepsilon(x)\overline{A}(1-x).
\end{align*}
\end{cor}
\section*{Acknowledgement}
This work  was  supported by the Initial Foundation for Scientific Research  of Northwest A\&F University (No. 2452015321).

\end{document}